\documentclass[12pt]{amsproc}
\usepackage[cp1251]{inputenc}
\usepackage[russian]{babel}
\usepackage{amssymb}
\numberwithin{equation}{section}

\def\al{\alpha}
\def\la{\lambda}

\begin{document}
УДК 517.984.54

\title[]{Обратные задачи для оператора\\
Штурма-Лиувилля с потенциалами из пространств
Соболева.  Равномерная устойчивость. }
\thanks{Работа поддержана Российским фондом фундаментальных
исследований, грант № 10-01-00423 }

\author[А.М.Савчук,A. A. Шкаликов] {А.М.Савчук, А.А.Шкаликов}


\maketitle
{\bf  Аннотация}. {\it В работе изучаются  две обратных
задачи для оператора Штурма-Лиувилля $Ly=-y'' +q(x)y$
на отрезке $[0,\pi]$. С первой из них при $\theta
\geqslant 0$ связано отображение  $F:\, W^{\theta}_2
\to l^{\,\theta}_B,
\ F(\sigma) =
\{s_k\}_1^\infty$, где $ W^\theta_2=
W^{\theta}_2[0,\pi]$ --- пространство Соболева,
$\sigma =\int q$ --- первообразная потенциала $q$,  а
$l^{\,\theta}_B$  --- специально построенное
конечномерное расширение весового пространства
$l^{\,\theta}_2$, куда помещаются регуляризованные
спектральные данные ${\bold s}=\{s_k\}_1^\infty$  для
задачи восстановления по двум спектрам. Подробно
изучаются  свойства отображения $F$. Основной
результат --- теорема о равномерной устойчивости. Он
состоит  в доказательстве равномерных оценок и снизу
и сверху нормы разности $\|\sigma -
\sigma_1\|_\theta$ через норму разности
регуляризованных спектральных данных $\|{\bold
s}-{\bold s}_1\|_\theta$, где норма берется в
$l^{\,\theta}_B.$ Аналогичный результат получен для
второй обратной задачи, которая связана с
восстановлением потенциала по спектральной функции
оператора $L$, порожденного краевыми условиями
Дирихле. Результат является новым и для классического
случая $q\in L_2$, который отвечает значению $\theta
=1$.}

\bigskip

В этой работе мы  изучим две классические обратные
задачи для оператора Штурма-Лиувилля
\begin{equation}
\label{1}
Ly =-y'' +q(x)\, y, \qquad  x\in [0,\pi],
\end{equation}
на конечном интервале. Первая задача связана с
восстановлением потенциала по двум спектрам этого
оператора, порожденного краевыми условиями Дирихле и
Дирихле-Неймана соответственно (мы называем ее
задачей Борга).  Вторая задача связана с
восстановлением потенциала по спектральной функции
этого оператора, порожденного краевыми условиями
Дирихле (далее такой оператор называем оператором
Дирихле). Давно известно решение этих задач для
вещественных потенциалов $q\in L_2$, в частности,
получена полная характеризация спектральных данных
для потенциалов $q$  из этого класса.  Наша цель ---
решить эти задачи для потенциалов $q$ из шкалы
соболевских пространств $ W^\alpha_2$ при всех
фиксированных $\alpha \geqslant -1$ (включая случай
$\alpha\in [-1,0)$,  когда потенциал является
сингулярной функцией-распределением.) Важную роль при
этом играют специальные гильбертовы пространства,
которые мы конструируем для решения указанных задач.
Эти пространства нужны для задания и изучения
отображений, которые мы связываем с рассматриваемыми
задачами, а также  для полного описания
(характеризации) спектральных данных, когда
первообразная потенциала $\sigma  = \int q(t)\, dt$
пробегает множество вещественных функций пространства
$W^{\alpha +1}_2.$


После решения обратных задач возникает важная задача об априорных оценках: насколько мало изменяется первообразная
потенциала $q$ в норме пространства $W^{\alpha +1}_2$ при малом изменении спектральных данных в норме соответствующего
гильбертова пространства, куда эти данные  помещены. Априорные оценки ранее были известны в классическом случае (при
$\alpha =0$). Но это были оценки локального типа, в которых постоянные вместе с радиусом окрестности, где оценки
действуют, зависели от потенциала $q$. Поэтому эффективность локальных оценок мала. \emph {Основная цель этой работы
--- получить равномерные двусторонние априорные
оценки не только для классического случая $\alpha =0$, но и при всех $\alpha > -1$}. Случай $\alpha = -1$  особый.
Развиваемый нами метод при $\alpha = -1 $ не работает.  Одновременно мы выясним, для каких спектральных данных
константы в априорных оценках могут "<портиться"> (т.е. становиться большими или малыми). Такая информация важна для
реализации конкретных вычислений при решении обратных задач. Мы покажем, что константы в априорных оценках
"<ухудшаются"> только по двум причинам: 1) увеличение нормы регуляризованных спектральных данных, т.е. большие
уклонения спектральных данных от нулевых значений (которые соответствуют нулевому потенциалу $q$); 2) уменьшение зазора
(расстояния) между парами соседних собственных значений или приближение к нулю одного из нормировочных чисел
(уменьшение числа $h$, которое фигурирует в определении множеств регуляризованных спектральных данных $\Omega_B^\theta
(h, r)$ и $\Omega_D^\theta (h, r)$, определяемых в параграфах 2 и 3 работы). Отметим еще раз, что \emph {эти оценки
являются новыми и для классического случая потенциалов $q\in L_2$,}  но метод доказательства оценок существенно
использует предварительные результаты, полученные при изучении обратных задач для потенциалов $q$  во всей шкале
соболевских пространств $W_2^\alpha$.

 История изучения обратных задач для оператора Штурма--Лиувилля
 ведет начало от работы Амбарцумяна \cite{Am}.
 Но результат этой работы  оказался не характерным для
 теории. Пионерскую роль сыграла  фундаментальная  работа Борга
\cite{Bo}, основной результат которой --- теорема единственности
для восстановления  потенциала по двум спектрам.
 Другую интерпретацию результатов Борга предложил
 Левинсон \cite{Lev}.  Тихонов \cite{Ti} показал,
 что потенциал  (при некоторых дополнительных условиях)
   восстанавливается единственным образом
 по функции Вейля--Титчмарша.  Марченко \cite{Ma1, Ma2}
  первым применил в исследовании  обратных задач
  оператор преобразования и доказал единственность
  решения обратной задачи по спектральной функции
  для операторов Штурма-Лиувилля, как на конечном интервале, так и на
  всей оси. Гельфанд  и Левитан \cite{GL} нашли необходимые и
  достаточные условия для восстановления потенциала
  по спектральной функции  и написали явные уравнения
  для решения задачи о восстановлении. Левитан
  \cite{Le1},  а также  Гасымов и Левитан
  \cite{GaLe},
   получили аналогичные результаты для задачи Борга
   о восстановлении потенциала по двум спектрам.
   Полное решение задачи Борга для потенциалов из
   $L_2$  получил Марченко \cite{Ma3}.
   Другие формулы для решения обратных задач
   предложил Крейн \cite{Kr1, Kr2}.
   В серии работ Трубовица  с соавторами был предложен
метод    для решения некоторых обратных   задач на
конечном интервале,
   использующий язык теории
   аналитических отображений. Детальное изложение
   имеется в книге Пошеля и Трубовица \cite
{PoTr}. Из последних работ, развивающих этот метод,
отметим работу Коротяева и Челкака \cite{KCh}.
  Для решения нелинейных
уравнений важную роль сыграла обратная задача по
данным рассеяния, изучение которой было проведено
Фаддеевым \cite{Fa1, Fa2}, Дейфтом и Трубовицем
\cite{DeTr}, Марченко \cite{Ma3}
(см. \cite{Ma3} для более полной информации). Большое
число работ посвящено изучению прямых и обратных
задач для операторов Штурма--Лиувилля в импедансной
форме. Отметим, что имеется связь между такими
операторами и обычными операторами Штурма-Лиувилля с
сингулярными потенциалами.
Работа Альбеверио, Гринива и Микитюка
\cite{AHM}
--- одна из последних на эту тему;  в ней имеются
многочисленные ссылки.

 В работе \cite{SS1}  авторы
предложили метод регуляризации для определения
оператора Штурма-Лиувилля с
потенциалами-распределениями $q\in W^{-1}_2$. Гринив
и Микитюк \cite{HM1, HM4} показали существование
оператора преобразования для уравнений с такими
потенциалами и дали решения классических обратных
задач для потенциалов $q\in W^{-1}_2$ (cм., в
частности,
\cite{HM2, HM3, HM5}). Марченко и Островский
 \cite{MaO}, \cite{Ma3} дали описание спектральных
 данных задачи Борга для потенциалов $q$  из пространств
 Соболева $W^{\alpha}_2 $
 при целых показателях гладкости $\alpha =
 0,1,2,\dots $. Аналогичные результаты для обратных задач по
спектральным функциям получили Фрайлинг и Юрко
\cite{FrYu}. Авторы \cite{SS3}   ввели шкалу
пространств $l^{\,\alpha +1}_B$ для изучения
спектральных данных задачи Борга и в терминах этих
пространств провели исследование этой задачи при всех
показателях гладкости $\alpha \geqslant -1$. В других
терминах и другим методом задачу Борга, а также
обратную задачу по спектральной функции исследовали
для показателей гладкости $\alpha\in
 [-1,0]$ Гринив и Микитюк \cite{HM6}.

 Различные априорные оценки локального характера для
обратных задач получали многие авторы.
Не вникая в детали, отметим, что в этом направлении
результаты  получили Марченко и Маслов
\cite{MaMa}, Рябушко \cite{Rja1, Rja2},
Хохштадт \cite{Ho},   Хальд
\cite{Ha},  Юрко \cite{Yu} (см.также \cite {FrYu}), Мизутани \cite{Mi},
Алексеев \cite{Al}, Мак Лафлин \cite{McL}, Хитрик
\cite{Hi}, Марлетта и Вайкард \cite{MaWe}, Маламуд
\cite{Ma}.

 Говоря об обратных задачах на конечном интервале,
необходимо упомянуть задачу о восстановлении
потенциала по двум спектрам периодической и
антипериодической задач. Естественно, она связана с
изучением оператора Хилла на всей прямой. С этой
задачей связано много интересных работ, и  изучена
она наиболее полно. Выделим важные результаты
Марченко и Островского \cite{MaO} и \cite{MaO2}.
Обратная задача для периодического случая является
единственной, для которой  получены равномерные
априорные оценки разности потенциалов через разность
спектральных данных (см. \cite{MaO2}).
  Из последних публикаций о периодической задаче отметим работу
Джакова и Митягина
\cite{DM1}, где наряду с новыми важными результатами имеется
подробная информация и библиография для случая
классических потенциалов, а также их работу
\cite{DM2}, где рассматриваются сингулярные
потенциалы.   Для более подробных сведений по
рассматриваемым здесь и другим обратным задачам мы
отсылаем читателя к монографиям Марченко
\cite{Ma3}, Левитана \cite{Le2}, Фрайлинга и Юрко
\cite{FrYu}, а также к обзорной работе Гештези
\cite{Ge}.

 Настоящая работа является продолжением серии работ
авторов \cite{SS3}--\cite{SS5}, посвященных решению
обратных задач с потенциалами из пространств
Соболева. В этих работах уже были сконструированы
пространства, в которые следует помещать
регуляризованные спектральные данные рассматриваемых
двух обратных задач, и изучены свойства отображений,
ставящих в соответствие первообразной потенциала
$\sigma =\int q(t)\, dt$ регуляризованные
спектральные данные. Ключевой является формулируемая
ниже в нужной форме Теорема 1.3 о слабой нелинейности
построенных отображений (эта теорема для разных задач
доказана в работах
\cite{SS4}-\cite{SS5}).
Как уже говорилось,  решение задачи Борга для
потенциалов $q\in W^\alpha_2$  во всей шкале $\alpha
\geqslant -1$  было дано авторами \cite{SS3}.
 Инъективность
доказывалась модификацией метода Борга, а для
описания образа и процедуры восстановления
развивались идеи работ Трубовица и соавторов. В этой
работе мы дополним исследования \cite{SS3} по задаче
Борга, в частности, докажем  локальную устойчивость
для всех индексов гладкости $\alpha \geqslant -1$. Мы
покажем также, что решение в пространствах Соболева
обратной задачи по спектральной функции оператора
Дирихле   для всех индексов $\alpha
\geqslant -1$ может быть проведено по такой же схеме,
как решение задачи Борга. Однако при реализации этой
схемы доказательства некоторых похожих утверждений
требуют новых подходов. \emph { Но главная цель
--- равномерные априорные оценки, которые мы получаем при
$\alpha >-1$.} Для их доказательства развивается
новый метод, основанный на теоремах о слабой
нелинейности построенных отображений.

Первый параграф работы носит вспомогательный
характер. Мы напоминаем основные определения и
конструкции пространств и формулируем в нужном виде
необходимые для дальнейшего результаты работ
\cite{SS3}-\cite{SS5}. Второй параграф основной. Здесь
приводятся дополнения к результатам
\cite{SS3} по задаче Борга, доказываются локальные и
равномерные априорные оценки для этой задачи. В
третьем параграфе все результаты для задачи Борга
переносятся на обратную задачу по спектральной
функции оператора Дирихле.

\bigskip

{\bf \section{ Определения пространств и нелинейных отображений, связанных с обратными задачами. Теоремы о свойствах
таких отображений.}}


\vspace{0.2cm}

 Сначала напомним, что определение
оператора Штурма-Лиувилля  с классическим потенциалом
$q\in L_1[0,\pi]$ можно распространить на случай
потенциалов-распределений из cоболевского
пространства $W^{-1}_2[0,\pi]$. Предположим, что
комплекснозначный потенциал $q$ принадлежит
соболевскому пространству $W^\alpha_2[0,\pi]$ при
некотором $\alpha
\geqslant -1$. Положим $\sigma(x) = \int q(x)\, dx$,
где первообразная понимается в смысле теории
распределений. Согласно
\cite{SS1} (см. также  \cite{SS2}, где даны альтернативные определения),
 определим
оператор Дирихле равенством
\begin{equation}
\label{2}
L_Dy=Ly=-(y^{[1]})'-\sigma(x)y^{[1]}-\sigma^2(x)y,\quad
y^{[1]}(x):=y'(x)-\sigma(x)y(x),
\end{equation}
взяв в качестве области определения
$$
\mathcal D(L_D)=\{y,\, y^{[1]}\in W_1^1[0,\pi]\ \,\vert\ \,  Ly\in L_2[0,\pi],\ y(0)=y(\pi)=0\}.
$$
Оператор Дирихле--Неймана определим аналогично:
$L_{DN} y= Ly$ на области
$$
\mathcal D(L_{DN})=\{y,\,y^{[1]}\in W_1^1[0,\pi]\ \,\vert\
\ \, Ly\in L_2[0,\pi],\ y(0)=y^{[1]}(\pi)=0\}.
$$

Для гладких функций $\sigma$ правые части в
\eqref{1} и \eqref{2} совпадают, и мы получаем
в первом случае  классический оператор
Штурма--Лиувилля с краевыми условиями Дирихле, а во
втором случае оператор с краевым условие $y(0) =0$  и
смешанным краевым условием $ y'(\pi) -h y(\pi) =0$,
где $h=\sigma(\pi)$. В первом случае оператор не
зависит от выбора константы в определении
первообразной $\sigma$  потенциала $q$, а во втором
случае зависит. Если константу выбрать так, чтобы
$\sigma(\pi) =0$, то мы получим классический оператор
Дирихле--Неймана.

Теперь определим спектральные данные для
рассматриваемых в работе задач. Обозначим через
$s(x,\lambda)$ единственное решение уравнения
$Ly-\lambda y =0$, удовлетворяющее условиям
 $s(0,\lambda) = 0$ и $s^{[1]}(0,
\lambda)= \sqrt\lambda$ (известно \cite{SS1}, что
такое решение существует и единственно). Очевидно,
что нули  $\{\lambda_k\}_1^\infty$ и
$\{\mu_k\}_1^\infty$ целых функций
$s(\pi,\lambda)/\sqrt\lambda$ и
$s^{[1]}(\pi,\lambda)/\sqrt\lambda$ являются
собственными значениями операторов $L_D$ и $L_{DN}$
соответственно. В случае вещественного потенциала $q$
все нули этих функций  являются простыми и
вещественными, и мы считаем их занумерованными так,
чтобы обе последовательности были строго
возрастающими. Для комплексных $q$ нумерацию можно
провести так, чтобы последовательности
$\{|\lambda_k|\}_1^\infty$ и $\{|\mu_k|\}_1^\infty$
не убывали.

 В задаче Борга потенциал восстанавливается по двум
спектрам
 $\{\lambda_k\}$ и $\{\mu_k\}$  операторов
 $L_D$  и $L_{DN}$.
Задание этих двух спектров эквивалентно заданию чисел
\begin{equation*}
s_{2k-1}=  \sqrt{\mu_k } -(k-1/2),\qquad s_{2k}
=\sqrt{\lambda_k } -k,\qquad k=1,2, \dots ,
\end{equation*}
т.е. последовательности
$\{s_k\}_1^\infty=\{s_k(B)\}_1^\infty$. Будем
говорить, что такая последовательность определяет
\emph{регуляризованные спектральные данные} задачи
Борга. Здесь  и далее мы подразумеваем,  что в
приведенных формулах ветвь квадратного корня выбрана
так, что значения аргумента $\sqrt\lambda$ лежат в
сегменте $(-\pi/2,\,\pi/2]$.

 Известно  \cite[Гл.3]{Le1}, что спектральная
функция оператора Дирихле однозначно
восстанавливается по его собственным значениям и так
называемым {\it нормировочным константам}, которые
определяются равенствами
\begin{equation*}
\alpha_k =\begin{cases}\int_0^\pi s^2(x,\lambda_k)\, dx,\quad\text{если}\ \la_k\ne0;\\\phantom{.}\\
\int_0^\pi\left(\tfrac{s(x,\lambda)}{\sqrt{\la}}\right)^2\, dx\Big\vert_{\la=\la_k},\quad\text{если}\ \la_k=0.
\end{cases}
\end{equation*}
Такое определение нормировочных чисел мы сохраним и для комплексных потенциалов. Последовательности
$\{\lambda_k\}_1^\infty \cup \{\alpha_k\}_1^\infty$ формируют спектральные данные оператора $L_D$. Задание этих данных
эквивалентно заданию чисел
\begin{equation}\label{L_D}
s_{2k} = \sqrt{\lambda_k } -k,\qquad s_{2k-1} =
\alpha_k -\pi/2,\qquad k=1,2, \dots,
\end{equation}
Будем говорить,что последовательность
$\{s_k\}_1^\infty =\{s_k(D)\}_1^\infty$ определяет
{\it регуляризованные спектральные
данные оператора $L_{D}$}.

 Имеем две задачи:  восстановить первообразную
потенциала $q$   по регуляризованным спектральным данным либо оператора $L_D$,  либо задачи Борга. Ясно, что в
сингулярном случае  восстановление функции $q$ невозможно  и работать надо с ее первообразной $\sigma =
\int q(x)\, dx$. При $q\in W^\alpha_2,
\alpha
\geqslant -1$  имеем $\sigma
\in W^\theta_2$, где $\theta =\alpha +1
\geqslant 0$. Случай классического потенциала $q\in L_2$
  соответствует показателю $\theta =1$. Нужно еще
  отметить, что переход  к восстановлению первообразной
меняет постановку задачи. Например, при
восстановлении  дифференцируемой функции $\sigma$ по
спектральным данным задачи Борга восстанавливается не
только потенциал $q = \sigma'$, но и постоянная
$h=\sigma(\pi)$
 в смешанном краевом условии.  Но по спектральным
 данным оператора $L_D$  функция $\sigma$
 восстанавливается только с точностью до постоянной.

Чтобы далее использовать язык теории отображений,
нужно понять, каким пространствам принадлежат
определенные выше  регуляризованные спектральные
данные, когда первообразная $\sigma$ пробегает
соболевское пространство $W^\theta_2, \theta\geqslant
0$.  Ясно, что эти пространства  разные для
рассматриваемых нами  двух задач. Однако различаются
они незначительно. Для  обеих задач эти пространства
можно выбрать, как конечномерные расширения обычных
весовых $l_2-$ пространств. Как расширять ---
становится ясным после анализа асимптотических формул
для собственных значений $\lambda_n, \
\,\mu_n$  и нормировочных констант
$\alpha_n$. Подробно это объяснено в
\cite{SS4}, \cite{SS5}.  Понять это можно также из
формулируемой ниже Теоремы 1.2  после интегрирования
по частям  формул, которыми определены операторы
$T_B$ и $T_D$.

Построим пространство для регуляризованных
спектральных данных задачи Борга. Обозначим через
$l^{\,\theta}_2 $ весовое $l_2$-пространство,
состоящее из последовательностей $\bold x=\{x_1, x_2,
\dots\}$,комплексных чисел,
 таких, что
$$
\|\bold x\|^2_\theta : =\sum_1^\infty |x_k|^2\, k^{2\theta} <\infty.
$$
Рассмотрим специальные последовательности
$$
{\bold e}_{2s-1} = \{k^{-(2s-1)}\}_{k=1}^\infty, \qquad
{\bold e}_{2s} = \{(-1)^k\, k^{-(2s-1)}\}_{k=1}^\infty,
\qquad s=1,2,\dots .
$$
Пусть $m=[\theta/2+3/4]$, где $[a]$ ---  целая часть
числа $a$. Положим
$$
 l_{B}^{\, \theta}=l_2^{\, \theta}\oplus
span \{{\bold e}_{ k}\}_{k=1}^{2m}.
$$
Здесь мы учли, что при $k\leqslant 2m$
последовательности ${\bold e}_k$  не принадлежат
пространству $l^{\,\theta}_2$, а при $k>2m$
принадлежат. Таким образом, \ $l_B^{\,\theta}$
состоит из элементов $\bold x +
\sum_{k=1}^{m} c_k \bold e_k$, где $\bold x \in
l_2^{\,\theta}$, а $\{c_k\}_1^{m}$ --- произвольные
комплексные числа. Скалярное произведение элементов
из $ l_B^\theta$ определяется формулой
$$
(\bold{x} +\sum_{k=1}^{m} c_k \bold e_k , \ \bold{y}
+\sum_{k=1}^{m} d_k \bold e_k ) = (\bold x ,
\bold y)_\theta + \sum_{k=1}^{m}
c_k\overline{d_k}.
$$

Построенное пространство  свяжем с регуляризованными
спектральными данными для оператора $L_B$. Хотя это
пространство определено как конечномерное расширение
весового пространства $l_2^{\,\theta}$, его элементы
удобнее записывать в форме обычных
последовательностей. Например, при $3/2\leqslant
\theta <5/2$  пространство $l_B^{\,\theta}$
состоит из последовательностей ${\bold
x}=\{x_k\}_1^\infty$ с координатами
$$
x_k =y_k +\alpha_1 k^{-1} +\alpha_2 (-1)^k k^{-1},
\quad
 \text{где}\ \ \{y_k\}_1^\infty \in
l_2^\theta, \ \ \alpha_1, \alpha_2 \in \mathbb C.
$$
Из такого представления легко следует, что
пространство $l_D^{\,\eta}$  компактно вложено в
пространство $l_D^{\,\theta}$  при $\eta > \theta$
(здесь мы принимаем во внимание компактность вложения
$l^{\,\eta}_2\hookrightarrow l_2^{\,\theta}$\ при
$\eta > \theta$).

Для построения пространства $l^{\,\theta}_D$
регуляризованных спектральных данных для оператора
Дирихле нужно вместо последовательностей  ${\bold
e}_k$ использовать последовательности
$$
 \widehat{\bold e}_{2s-1} =
\{\, 0,\  2^{-(2s-1)},\, 0,\ 4^{-(2s-1)},\, 0,\
6^{-(2s-1)},\ldots\},
 \quad \widehat{\bold e}_{2s}=
\{2^{-(2s)},\, 0,\  4^{-(2s)},\, 0,\  6^{-(2s)},\ldots\}.
$$
Пространство  $l_D^{\, \theta}$  определим равенством
$l_D^{\, \theta} = l_2^{\,
\theta}\oplus span \{\widehat{\bold e}_{ k}\}_{k=1}^m,$
где число $m$ однозначно определено условием $m-1/2
\leqslant \theta < m+1/2$.
 Отметим, что в работе
\cite{SS5}
конструкция пространства для регуляризованных
спектральных данных оператора $L^\theta_D$
проводилась в пространстве двусторонних
последовательностей. Здесь мы реализовали
эквивалентную конструкцию в пространстве
односторонних последовательностей, чтобы оба
пространства выглядели единообразно.

Определим следующие нелинейные операторы:
\begin{equation}\label{defi}
F_B (\sigma) = \{s_k(B)\}_1^\infty ,
\quad F_D (\sigma) =
\{s_k(D)\}_1^\infty.
\end{equation}
Из результатов работ
\cite{SS2} и
\cite{HM1}
следует, что последовательности, образованные из
регуляризованных спектральных данных  в правых частях
равенств \eqref{defi}, являются последовательностями
из $l_2$ для любой первообразной  $\sigma =\int
q(x)\, dx\in L_2(0,\pi)$. Поэтому все выписанные в
\eqref{defi} операторы корректно определены как
операторы из $L_2$
 в $l_2$. Более того,  согласно результатам из \cite{SS3} и \cite
 {SS5}, образы сужений этих операторов на соболевские
 пространства $W^\theta_2,\ \, \theta >0,$
 лежат в пространствах $l_B^{\,\theta}$\  и \  $l_D^{\,\theta}$
 соответственно. Именно для этой цели мы проводили
 расширения пространств  $l_2^{\,\theta}$. Без
 присоединения к  $l_2^{\,\theta}$  специальных последовательностей
 соответствующий результат неверен.

Далее  будут использоваться результаты работ
\cite{SS3}-\cite{SS5},  которые приведем
 в нужном нам виде.

{\bf Теорема 1.1} {\sl При любом фиксированном $\theta \geqslant 0$
нелинейные операторы $F_B $ и $F_D$ корректно
определены как операторы из пространства $W^\theta_2$
в  $l_B^{\,\theta}$  и
\ $l_D^\theta$ соответственно. Эти операторы дифференцируемы по
Фреше в каждой точке (функции) $\sigma$ при условии,
что эта функция вещественнозначна  и все собственные
значения $\lambda_k(\sigma),\  \mu_k(\sigma)$ не
обращаются в нуль (для отображения $F_D$ достаточно,
чтобы не обращались в ноль только
$\lambda_k(\sigma)$). В частности, эти операторы
дифференцируемы по Фреше в точке $\sigma=0$, причем
производные по Фреше в этой точке суть линейные
операторы $T_B$ и $  T_D$, которые определяются
формулами}
\begin{equation*}
(T_B\,\sigma)_k=-\frac1{\pi}\int\limits_0^\pi\sigma(t)
\sin(kt)\, dt,
\quad k=1,2,\dots.
\end{equation*}

\begin{equation*}
\left\{ \begin{array}{l} (T_D\,\sigma)_{2k-1}
=-\int\limits_0^\pi(\pi-t)
\sigma(t)\cos(2kt)\, dt ,
\quad k=1,2,\dots,\\
(T_D\,\sigma)_{2k}=-\frac1{\pi}\int\limits_0^\pi\sigma(t)
\sin(2kt)\, dt,\quad k=1,2,\dots .
\end{array}\right.
\end{equation*}

{\bf  Доказательство}  этого утверждения для
 оператора $F_B$ получается из Предложения 1 и Теоремы
6.1 работы
\cite{SS4}, а для второго  оператора --- из Предложения
1 и Теоремы 4.2 работы  \cite {SS5}.\quad$\Box$

{\bf Теорема 1.2}.{\sl  Пространства $l_B^{\,\theta}$\  и \
 $l_D^\theta$ образуют шкалу
компактно вложенных друг в друга пространств,
замкнутых относительно интерполяции, т.е. $[l^0,\,
l^{\,\theta}]_\tau =l^{\,\theta\tau}$ при всех
$\theta
\geqslant 0,
\tau \in [0,1]$ (здесь для краткости опущены нижние  индексы
$B$  или $D$). При любом $\theta
\geqslant 0$ оператор $T_B$ изоморфно отображает
пространство $W_2^\theta$ на \ $l_B^{\,\theta}$.
Оператор
  $T_D$\ изоморфно отображает
пространство \ $W_2^\theta \ominus\{1\}$  на
\ $l_{D}^{\,\theta}$.}

{\bf Доказательство.}
Первое утверждение этой теоремы для пространства
$l^\theta_B$  доказано в Предложении 4 работы
\cite{SS4}. Доказательство для  пространства $l^\theta_D$
проходит без изменений. Второе утверждение для
оператора $T_B$   доказано в Лемме 1 работы
\cite{SS4},  а для оператора $T_D$ в Предложении 3 работы
\cite{SS5}.\quad$\Box$

Следующая теорема является наиболее существенным звеном в доказательстве основных результатов этой работы. В частности,
она говорит, что рассматриваемые отображения $F_B$ и $F_D$  являются слабо нелинейными, т.е. компактными возмущениями
линейных отображений. Важна также точная зависимость от $\theta$ показателя $\tau =\tau (\theta)$, который
характеризует "<качество"> компактности.

 {\bf Теорема 1.3.} {\sl При любом фиксированном $\theta \geqslant 0$
 оператор \ $F_B$\  отображает пространство $W^\theta_2$\
 в \  $l_B^{\,\theta}$\  и допускает представление вида
 $$
 F_B(\sigma) = T_B\ \sigma + \Phi_B(\sigma).
$$
Здесь  \ $T_B$ \ --- линейный оператор, определенный
в Теореме 1.1, а $\Phi_B$\ отображает пространство \
$W_2^\theta$
\ в $l^{\,\tau}_B$, \ где
$$
\tau=\begin{cases}2\theta,\quad& \text{если}\
0\leqslant\theta\leqslant1,\\ \theta+1,\quad& \text{если}\
1\leqslant\theta<\infty.\end{cases}
$$
Кроме того, отображение  \ $\Phi_B:\,W_2^\theta\to
l_B^\tau$\ является ограниченным в любом шаре, т.e.
\begin{equation*}
\|\Phi(\sigma)_B\|_\tau\leq C(R),\quad\text{если}\ \|\sigma\|_\theta\le R,
\end{equation*}
где постоянная \ $C$\ зависит только от радиуса шара
$R$.
\ Аналогичное утверждение справедливо для оператора
\ $F_D$.\ А именно,
$$
  F_D(\sigma) = T_D\ \sigma + \Phi_D(\sigma)
$$
и отображение
\ $\Phi_D:\,W_2^\theta\ominus \{1\} \to l_D^{\,\tau}$\
обладает тем же свойством, что и $\Phi_B$.}

{\bf Доказательство} этой теоремы для оператора $F_B$
проведено в работе \cite{SS4},  а для второго
оператора --- в работе \cite{SS5}. В случае $\theta
>0$ компактность нелинейных слагаемых в представлениях
операторов $F_B $ и $F_D$ вытекает из компактности
вложений \ $l^{\,\eta}
\hookrightarrow
l^{\,\theta}$\ при условии $\eta
>\theta$\ (здесь мы опускаем для краткости индекс
$B,$ \ или \ $D$). Случай $\theta =0$ особый. При $\theta =0 $  из  сформулированной теоремы не вытекает компактность
нелинейных слагаемых.\quad$\Box$

\vspace{0.2cm}

\section{Задача Борга. Характеризация  спектральных данных
для  первообразных $\sigma$ вещественных потенциалов
$q\in W^\alpha_2$. Равномерные априорные оценки.}

\medskip
В этом   и следующем параграфе мы используем
следующие обозначения. Через $W^\theta_{2,\mathbb R}$
обозначаем множество вещественных функций в
пространстве $W^\theta_2$, через $\mathcal B^{\,
\theta}_{\mathbb R} (R)$ ---  замкнутый шар радиуса
$R$ в $W^\theta_{2,\mathbb R}$, через
 $\Gamma^{\, \theta}_B$ ---
 множество всех функций в
$W^\theta_{2,\mathbb R}$, для которых
$\mu_1(\sigma)\geqslant 1/4$,  и через $\mathcal
B_\Gamma^{\,\theta} (R)$
--- пересечение множества $\Gamma^{\, \theta}_B$ и шара
$\mathcal B^{\, \theta}_{\mathbb R} (R)$.  Здесь
$\mu_1(\sigma)$
--- первое собственное значение оператора $L_{DN}$.
Число $1/4$ взято для определенности и простоты,
вместо $1/4$ может фигурировать любое число $\eta>0$,
но тогда в \eqref{repll}  и
\eqref{replh} нужно писать $s_1 \geqslant \sqrt\eta
-1/2$.

Известно, что для вещественных потенциалов спектры
$\{\la_k\}$ и $\{\mu_k\}$ операторов $L_D$ и $L_{DN}$
удовлетворяют условию перемежаемости
\begin{equation}\label{repl}
\mu_1<\la_1<\mu_2<\la_2<\dots <\mu_n<\la_n<\mu_{n+1}<
\dots .
\end{equation}
Для классических потенциалов этот факт известен давно
(см., например, \cite{Ma2}), а для сингулярных
потенциалов-распределений он доказан в \cite{HM3} и в
\cite{SS3}.
Заметим, что для положительных $\lambda_k$ и $\mu_k$
неравенства \eqref{repl} эквивалентны неравенствам
для корней из этих чисел. Поэтому условия
\eqref{repl} вместе с условием $\mu_1 \geqslant 1/4$
(т.е. условием $\sigma
\in \Gamma^0_B$)  эквивалентны неравенствам
\begin{equation}\label{repll}
s_1 \geqslant 0, \qquad s_k-s_{k+1}<\frac 12, \qquad
k=1,2,\dots ,
\end{equation}
где $\{s_k\} =\{s_k(B)\}$ ---  регуляризованные
спектральные данные для задачи Борга.
Последовательность $\{s_k\}_1^\infty$ принадлежит
$l_2$, поэтому  для любой фиксированной вещественной
функции $\sigma\in L_2$ найдется число $h = h
(\sigma)>0$, такое, что
\begin{equation}\label{replh}
s_1 \geqslant 0, \qquad s_k-s_{k+1}\leqslant \frac
12-h,
\qquad k=1,2,\dots .
\end{equation}

 Фиксируем произвольные числа $r>0$ и
$h \in (0,\, 1/2)$. Обозначим через
$\Omega_B^{\,\theta} (r,h)$ совокупность вещественных
последовательностей $\{s_k\}_1^\infty$, для которых
выполняются неравенства
\eqref{replh}  и которые лежат в замкнутом шаре радиуса $r$
пространства $l_B^{\,\theta}$,  т.е.
$\|\{s_k\}\|_\theta \leqslant r$. Через
$\Omega_B^{\,\theta}$ обозначим множество  всех
вещественных последовательностей $\{s_k\}_1^\infty
\in l_B^{\,\theta}$, для которых справедливы неравенства
\eqref{repll}.

Напомним, что с задачей Борга мы связали оператор
$$
F_B: W^\theta_2 \to l^{\,\theta}_B, \qquad
F_B(\sigma) =
\{s_k\}_1^\infty,
$$
где $\{s_k\}^\infty_1$ ---  регуляризованные
спектральные данные задачи Борга. Из сказанного выше
 и Теоремы 1.3 следует, что $F_B$  отображает $\Gamma^{\,\theta}_B$ в
$\Omega_B^{\,\theta}$.

Далее в этом параграфе, там, где это удобно,  {\it мы
будем опускать индекс $B$}, так как будем работать
только с задачей Борга. В частности, операторы  $F_B,
T_B$ и $\Phi_B$ из Теоремы 1.3  обозначаем через $F,
T$ и $\Phi$ соответственно. Всюду вместо
$\Gamma^{\,\theta}_B,\ \, \Omega^{\,\theta}_B$ и
$\Omega^{\,\theta}_B(r,h)$ пишем $\Gamma^{\,\theta},\
\,\Omega^{\,\theta}$  и $\Omega^{\,\theta}(r,h)$.
Однако обозначение $l^{\,\theta}_B$ для пространств
регуляризованных спектральных данных  сохраняем
прежним.

{\bf Теорема 2.1} {\sl
При любом фиксированном $\theta
\geqslant 0$\ отображение \ $F:
\Gamma^{\,\theta} \to \Omega^{\,\theta}$\ есть биекция.
 }

{\bf Доказательство.}  Инъективность отображения
$F:\Gamma^{\,\theta} \to \Omega^{\,\theta}$\ доказана
в Лемме 6  работы
\cite{SS3}. В Лемме 5 этой же работы  проведено доказательство
 и сюръективности
этого отображения, но оно нуждается в дополнении в
случае  $\theta \geqslant 1/2$.  При $\theta <1/2$
пространство $l^{\,\theta}_B $ совпадает с
$l^{\,\theta}_2 $, а при  $\theta
\geqslant 1/2$ содержит еще состоящее из специальных
последовательностей подпространство $\mathcal L^{2m}$
размерности $2m$, где $m=[\theta/2 +3/4]$.
Познакомившись с доказательством Леммы 5 из
\cite{SS3}, приходим к выводу, что для полного его завершения
нужно уметь восстанавливать функцию $\sigma$ (или
доказывать ее существование), если варьируются только
координаты подпространства $\mathcal L^{2m}$, а все
координаты в $l^{\,\theta}_2$  остаются неизменными.
Авторы не видят простого прямого решения этой задачи,
без использования трудоемких теорем.  Здесь мы
приведем доказательство сюръективности с
использованием Теоремы 1.3, основываясь  на том, что
при $\theta\in [0,\, 1/2)$ это свойство уже доказано.

Нам известно, что отображение $F:
\Gamma^{\,\theta} \to \Omega^{\,\theta}$\
  сюръективно при $\theta \leqslant 1/4$. Покажем,
что оно сюръективно при любом $\theta\in (1/4,1/2]$. Возьмем произвольный элемент $\bold y\in
\Omega^{\,\theta}\subset l^{\,\theta}_B, \ \,
\theta \in(1/4,1/2]$. Поскольку при $\theta =1/4$
рассматриваемое отображение есть биекция, найдется
единственная функция  $ \sigma\in
\Gamma^{1/4}$, такая, что $F \sigma
=\bold y$ (здесь мы учитываем вложение
$l_B^{\,\theta}\hookrightarrow l_B^{1/4}$ ). В силу
Теоремы 1.3 имеем $T \sigma = -\Phi \sigma + \bold y
\in l_B^{\,\theta}$,  так как $\bold y \in l^{\, \theta}_B,$
а из условия $\sigma \in W^{1/4}_2$ следует, что
$\Phi\sigma\in l^{1/2}_B\hookleftarrow
l^{\,\theta}_B$. Но в силу Теоремы 1.2 линейный
оператор $T:\, W^\theta_2\to l^{\,\theta}_B$ есть
изоморфизм. Следовательно, $\sigma \in
W^\theta_{2,\mathbb R}$, а потому с учетом включения
$
\bold y\in
\Omega^{\,\theta}$ имеем $\sigma \in \Gamma^{\,\theta}$.
 Тем самым, мы доказали, что отображение сюръективно при
$\theta
\in (1/4, 1/2]$. Теперь, зная, что $F:
\Gamma^{\,\theta} \to \Omega^{\,\theta}$\
 сюръективно при $\theta\in[0,1/2]$, с помощью
такого же приема покажем сюръективность при $\theta\in (1/2,\, 1].$ Повторив этот же прием, с помощью Теоремы 1.3
покажем сюръективность при $\theta\in (1,2]$. На $k+1$-м шаге получим сюръективность при $\theta\in (k-1,k]$. Здесь
число $k$ произвольно, поэтому утверждение справедливо при всех $\theta \geqslant 0$. Теорема доказана.\quad$\Box$

Обозначим через $\widehat \Omega^{\,\theta}_B$
множество последовательностей $\{s_k\}_1^\infty \in
l^{\,\theta}_B$, для которых числа
$$
\mu_k =(s_{2k-1} +k-1/2)^2, \qquad
\lambda_k =(s_{2k-1} +k-1/2)^2
$$
вещественны и подчинены условиям \eqref{repl}.

Заметим,  что если   к функции $\sigma$, которой
определяются операторы $L_D$ и $L_{DN}$, добавить
функцию $c(x-\pi)$, то эти операторы перейдут в $L_D
+c$ и $L_{DN} +c$ соответственно, т.е. их спектры
сдвинутся на $c$. Положим
\begin{equation}\label{shift}
s_{2k-1}(c) = \sqrt{\mu_k +c} -(k-1/2),\qquad
s_{2k}(c) = \sqrt{\lambda_k +c} -k
\end{equation}
Поскольку $c(x-\pi) \in W^\theta_2$ при всех
$\theta\geqslant 0$,  то $\{s_k(c)\}^\infty_1 \in
l^{\,\theta}_B$, если и только если
$\{s_k(0)\}^\infty_1 \in l^{\,\theta}_B$.
Следовательно, $\{s_k\}^\infty_1 \in \widehat
\Omega^\theta$ если и только если найдется $c\geqslant
0$,  такое, что  $\{s_k(c)\}^\infty_1 \in
\Omega^\theta$. Из сделанных замечаний следует

{\bf Теорема  2.2.} {\sl Отображение $F:
W^\theta_{2,\mathbb R} \to \widehat\Omega^{\,\theta}$
есть биекция. Последовательности чисел
$\{\mu_k\}^\infty_1$ и $\{\lambda_k\}^\infty_1$
являются спектрами операторов $L_D$ и $L_{DN}$, если
и только если они удовлетворяют условиям
перемежаемости
\eqref{repl} и $\{s_k\}^\infty_1 \in
\l^{\,\theta}_B$. }

Отметим, что  при натуральных $\theta =1,2,\dots$
Марченко и Островский
\cite{MaO}, \cite{Ma3} провели характеризацию
спектральных  данных для задачи Борга  в другой
форме, без  использования  пространств
$l^{\,\theta}_B $. Можно показать, что для таких
значений $\theta$  их результат с учетом теоремы
единственности Борга эквивалентен сформулированной
теореме.

 Далее существенно будут использоваться
аналитические свойства отображения $F$.  Мы предполагаем, что читатель знаком с определением производных по Фреше и
Гато для отображения $F: U\to H$, где $U$ --- открытое множество в $E$, а $E$  и $H$  --- сепарабельные гильбертовы
пространства. Для комплексных гильбертовых пространств производная по Фреше естественно определяется в комплексном
смысле. Отображение $F: U\to H$  называется аналитическим, если существует комплексная  производная по Фреше в каждой
точке $x\in U$. Производную по Фреше в точке $x$ далее обозначаем через $F'(x)$. Естественным образом определяется
понятие вещественного аналитического отображения, см., например,
\cite{PoTr}. Отображение $F: U\to H$ называется слабо
аналитическим, если в комплексном смысле
дифференцируемы по Гато координатные функции $(F(x),
e_k)$, где $\{e_k\}_1^\infty$  --- ортонормированный
базис пространства $H$. Известен результат
\cite{PoTr}, который значительно упрощает проверку
аналитичности отображения.

{\bf Предложение 2.3}. {\sl Если $F: U\to H$ --- слабо
аналитическое отображение и локально ограничено в
каждой точке $x\in U$,  то $F$ --- аналитическое
отображение.}

Далее мы будем работать с отображениями замкнутых
множеств. Чтобы не делать дополнительных объяснений,
всюду считаем, что отображение $F:\,  D\to H$
аналитично на $D$,  если найдется открытое множество
$U$, такое, что $U\supset D$ и $F:\, U\to H$
аналитично.

{\it Доказательство.} Утверждения этой теоремы
 доказаны в параграфе 5 работы
\cite{SS4}. Доказательства  основаны на Теореме 1.3 и
Предложении 2.2, если предварительно вычислить
производные координат.
 Здесь важно, что знаменатели в формуле
\eqref{F'}  в случае вещественной функции $\sigma$ не
обращаются в ноль. Согласно \cite{SS1} собственные функции непрерывно зависят от первообразной потенциала $\sigma$, а
потому числа $(y^2_k(x), 1)$ не обращаются в ноль в некоторой комплексной окрестности (нужно еще учесть асимптотики
$y_k$ при $k\to\infty)$). Теорема остается справедливой, если вместо условия $\sigma\in\Gamma^{\,\theta}$ потребовать,
чтобы $\sigma$ была вещественной и среди чисел $\{\rho_k\} $ нет равных нулю. Однако в этом случае вместо вещественной
аналитичности будет обычная.

{\bf Теорема 2.4.} {\sl Пусть $\theta\geqslant0$  и
$\sigma\in\Gamma^{\,\theta}$. Тогда найдется
комплексная окрестность $U\in W_2^\theta$ точки
$\sigma$, такая, что отображение $F: U \to
l_B^\theta$ дифференцируемо в комплексном смысле во
всех точках этой окрестности.
 Таким образом, отображение
 $F: \Gamma^{\,\theta}\to l_B^{\,\theta}$
является вещественно аналитическим. Этим же свойством обладает отображение $\Phi =F-T: \Gamma^{\,\theta}
\to l^{\,\tau}_B $, где $T=T_B$ и $\tau$ определены в
Теореме 1.3. Производная в точке $\sigma\in \Gamma$
определяется равенством
\begin{equation} \label{F'}
[F'(\sigma)] f = \left\{ -
\frac{\left(y'_k(x)y_k(x),
\overline{f(x)}\right)}{\rho_k\, (y_k^2(x), 1)}\right\}_{k=1}^\infty.
\end{equation}
Здесь
 $\rho_{2n-1} =
\sqrt{\mu_n },\ \,\rho_{2n} =
\sqrt{\la_n},\ \, y_{2n-1}(x)$ --- собственные функции оператора
$L_{DN},\ \,  y_{2n}$ --- собственные функции
оператора $L_{D}$, а  $f\in W^\theta_2$ --- функция,
на которую действует оператор $ F'(\sigma)$.}

{\bf Доказательство.} Утверждения этой теоремы
 доказаны в параграфе 5 работы
\cite{SS4}. Доказательства  основаны на Теореме 1.3 и
Предложении 2.2, если предварительно вычислить
производные координат.
 Здесь важно, что знаменатели в формуле
\eqref{F'}  в случае вещественной функции $\sigma$ не
обращаются в ноль. Согласно \cite{SS1}, собственные функции непрерывно зависят от первообразной потенциала $\sigma$, а
потому числа $(y^2_k(x), 1)$ не обращаются в нуль в некоторой комплексной окрестности (нужно еще учесть асимптотики
функций $y_k$ при $k\to\infty)$). Теорема остается справедливой, если вместо условия $\sigma\in\Gamma^{\,\theta}$
потребовать, чтобы $\sigma$ была вещественной и среди чисел $\{\rho_k\} $ не было равных нулю. Однако в этом случае
вместо вещественной аналитичности будет обычная.\quad$\Box$

{\bf Лемма 2.5}. {\sl Пусть функции $y_k(x)$,
участвующие в Теореме 2.4, нормированы условиями
$y^{[1]}_k(0) =1$. Тогда   система функций
\begin{equation}\label{phi}
\varphi_k(x) = \frac 2\pi y'_k(x)y_k(x), \qquad k=1,2, \dots,
\end{equation}
 является базисом Рисса в пространстве $L_2(0,
\pi)$.
 Биортогональная система к $\{\varphi_k(x)\}_1^\infty$
имеет вид
\begin{equation}\label{psi}
\psi_k(x)=\pi \rho^{1/2}_k y_k(x)w_k(x),
\end{equation}
где при $k=2n$ функция  $w_k$
--- решение уравнения $-y''+\sigma'y=\la_n y$ с
начальными условиями
$$
w_k^{[1]}(\pi)=0, \quad w_k(\pi)=\left(\int_0^\pi
y_k^2(x)dx\cdot y_k^{[1]}(\pi)\right)^{-1},
$$
а при $k=2n-1$ функция $w_k$ есть решение уравнения
$-y''+\sigma'y=\mu_n y$ с начальными условиями
$$
w_k(\pi)=0, \quad w_k^{[1]}(\pi)=-\left(\int_0^\pi
y_k^2(x)dx\cdot y_k(\pi)\right)^{-1}.
$$  }

{\bf Доказательство. }
 Первое утверждение
теоремы о базисности Рисса системы
$\{\varphi_k(x)\}_1^\infty$ доказано в Лемме 6 работы
\cite {SS3}. Там же доказаны соотношения
$(\varphi_k(x), \psi_m(x))=0 $ при $k\ne m$.
Доказательство равенств $(\varphi_k(x),
\psi_k(x)) =1$ проводится прямыми вычислениями,
которые мы здесь  опускаем, так как далее конкретный вид функций $\varphi_k$  и  $\psi_k$   не используется.\quad$\Box$

{\bf Теорема 2.6.} {\sl Пусть $\theta\geqslant 0$.
Для каждой точки $\bold y_0 \in
\Omega^{\,\theta} =F(\Gamma^{\,\theta})$
существует ее комплексная окрестность $U({\bold y_0})$, в которой определено обратное отображение $F^{-1}(\bold y)$ и в
которой это отображение имеет комплексную  производную  по Фреше. Эта производная имеет вид
\begin{equation}\label{F-1}
\left( F^{-1}\right)'(\bold y)= (F')^{-1}(\bold y) =
\sum^\infty_{k=1} s_k\widetilde\psi_k(x), \qquad \bold y= (s_1, s_2,
\dots ).
\end{equation}
Здесь $\widetilde\psi_k(x) =
\gamma_k \psi_k(x)$, где  $\{\psi_k(x)\}_1^\infty$
---  биортогональная система из Леммы 2.6, а $\gamma_k=
\rho_k\int_0^1 y^2_k(x)\, dx$.
 }

{\bf Доказательство.}
Пусть сначала $\theta >0$. Имеем
$$
F'(\sigma_0) = T + \Phi'(\sigma_0), \qquad
\bold y_0 = F( \sigma_0).
$$
 Согласно Теореме 1.2, оператор $T: W^\theta_2
\to l_B^{\,\theta}$
--- изоморфизм, а в силу Теоремы 2.4 оператор
$\Phi'(\sigma_0): W_2^\theta \to l_B^{\,\tau}$
ограничен, а потому оператор $\Phi'(\sigma_0):
W_2^\theta \to l_B^{\,\theta}$ компактен.
Следовательно, $F'(\sigma_0)$ --- фредгольмов
оператор, а потому он обратим, если его ядро нулевое.
 Из формул
\eqref{F'} и полноты системы \eqref{phi} в
пространстве $L_2$ следует, что равенство
$F'(\sigma_0)f =0$ при $f\in  L_2$ влечет $f=0$. Тем
более это так, если $f\in W^\theta_2$ при $\theta
>0$.
Формула \eqref{F-1} теперь получается
непосредственной проверкой. Достаточно проверить
равенство
$$
F'(\sigma_0)\left( F^{-1}\right)'(\bold y_0) =\bold
y_0.
$$
Оно сразу следует из \eqref{F'} и \eqref{F-1} c
учетом взаимной биортогональности систем
$\{\gamma_k^{-1}\varphi_k\}_1^\infty$ и
$\{\gamma_k\psi_k\}_1^\infty$.

Пусть теперь  $\theta =0$. Из асимптотических формул
для собственных значений $\rho_k^2$ и собственных
функций $y_k$, полученных в
\cite[Теоремы 2.6 и 2.7]{SS2},
 сразу следует, что $\gamma_k \asymp 1$, если
функции $y_k$ нормированы условием $ y_k^{[1]}(0)=1$.
Поэтому из Леммы 2.5 вытекает, что система
$\{\widetilde\psi_k\}_1^\infty$
--- базис Рисса. Тогда ограниченность оператора
$\left(F'\right)^{-1}(\bold y_0)$,  определенного формулой \eqref{F-1}, следует из  определения базиса Рисса.
Существование обратного оператора при любом $\theta\geqslant 0$ в малой комплексной окрестности точки $\bold y_0$  и
его комплексная дифференцируемость  следует из теоремы об обратном отображении. Теорема доказана.\quad$\Box$

 Отметим, что из теоремы 2.6  сразу получаются
локальные оценки разности потенциалов через разность спектральных данных и наоборот. Как отмечено во введении, для
классического случая $\theta =1\ (q\in L_2)$ имеется много работ на эту тему, выполненными различными методами и в
разной форме. Однако изучались отображения $q
\to \{\text{спектральные данные}\}$, мы же изучаем
отображение $\int q(t)\, dt = \sigma
\to \{\text{спектральные данные}\}$,
поэтому возникающие у нас системы и формулы имеют
другой вид.

Далее мы покажем, что при $\theta >0$  с помощью
Теоремы 1.3 можно получить существенно более сильный
результат, избегая технической работы с системами
функций.

{\bf Лемма 2.7.} {\sl Фиксируем $\theta >0$.
Пусть $R$ произвольное положительное число и
$\mathcal B^{\,\theta}_\Gamma (R) = \Gamma
\cap\mathcal B^{\,\theta}_\mathbb R (R)$.
 Тогда найдутся положительные числа $r=r(R), h=h(R)$,
такие, что
$$
F (\mathcal B^{\,\theta}_\Gamma (R) ) \subset
\Omega^\theta (r,h).
$$
}

{\bf Доказательство.}  Если    $\|\sigma\|_\theta \leqslant R,\
\sigma\in\Gamma^{\,\theta}$, то из Теоремы 1.3
следует, что $F\sigma = \bold y \in \Omega^{\,\theta}
(r) $, где $r=r(R)$ зависит  от $R$,  но не от
$\sigma$. Остается показать, что для всех элементов
$\bold y = F\sigma,\ \, \sigma
\in \mathcal B^{\,\theta}_\Gamma (R)$, выполняются
неравенства \eqref{replh} при некотором $h=h(R)>0$,
зависящем  от $R$,  но не от $\sigma$.


Заметим, что найдется число $N=N(\theta,r)$, такое,
что для всех $\bold y =(s_1, s_2, \dots )
\ \in
\Omega^{\,\theta} (r)$ и всех  $k\geqslant N$
выполняются неравенства $s_{k+1}-s_k \leqslant 1/4$
(здесь вместо $1/4$  можно взять любое число
$\varepsilon>0$).
 Это утверждение сразу следует из
определения нормы в $l_B^{\,\theta}$ при $\theta>0$
(см. подробнее \S \ 5 работы \cite{SS4}; при $\theta
=0$ это утверждение не справедливо). Теперь,
допустим, что утверждение теоремы неверно и найдутся
элементы ${\bold y}^n= F \sigma_n,
\ \sigma_n\in\mathcal B^{\,\theta}_\Gamma (R) $,
такие, что $s^n_k- s^n_{k+1}\to 1/2$  при
$n\to\infty$  и некотором фиксированном $1\leqslant k
< N$ (здесь $s^n_k $ --- координаты элементов $\bold
y^n$.)  Шар в  пространстве $W^\theta_2$  слабо
компактен, поэтому из последовательности функций
$\sigma_n$ можно выделить слабо сходящуюся
подпоследовательность. Не ограничивая общности,
считаем, что сама эта последовательность слабо
сходится к функции $\sigma\in W^\theta_{2,\,\mathbb
R}.$
 Так как пространство $W^\theta_2$ компактно вложено
в $L_2$, то последовательность $\sigma_n$ сильно
сходится к $\sigma$ в норме $L_2$. Пусть индекс $k$,
при котором $s^n_k - s^n_{k+1}\to 1/2$, является,
например, четным, $k=2p$. Тогда $\la_p(\sigma_n) -
\mu_{p+1}(\sigma_n) \to 0$. Согласно Теореме 2 работы
\cite{SS1}, сходимость функций  $\sigma_n$ в $L_2$
влечет за собой сходимость собственных значений, т.е.
$\la_p(\sigma_n)\to\la_p(\sigma),
\mu_{p+1}(\sigma_n)\to\mu_{p+1}(\sigma)$. Поэтому $s^n_k -
s^n_{k+1}\to 1/2$ влечет
$\la_p(\sigma)=\mu_{p+1}(\sigma)$, что невозможно в
силу условия перемежаемости \eqref{repl}. Лемма
доказана.

{\bf Лемма 2.8.} {\sl Пусть $\theta >0$. Справедливо обратное утверждение к Лемме~2.7: для любых чисел $r$ и $h$
найдется число $R>0$,  такое, что
$$
F^{-1}(\,\Omega^{\,\theta}(r,h)) \subset \mathcal
B^{\,\theta}_\Gamma (R).
$$
Справедливо представление
$$
F^{-1} = T^{-1} +\Psi, \qquad \Psi: \Omega^{\,\theta}
\to W^\tau_2,
$$
где число $\tau$  определено в Теореме~1.3.
Отображение $\Psi: \Omega^\theta \to W^\tau_2,$
аналитично, причем
\begin{equation}\label{Psi}
\|\Psi\bold y\|_\tau \leqslant C\|\bold y\|_\theta
\qquad\text{для всех} \ \, \bold y \in\Omega^{\,\theta}(r,h),
\end{equation}
где постоянная $C$  зависит только от $r$  и $h$. }

{\bf  Доказательство.}
Если первое утверждение леммы неверно,  то найдутся
элементы $\bold y^n \in
\Omega^{\,\theta}(r,h),$ такие, что $F^{-1}
\bold y^n = \sigma_n, \ \|\sigma_n\|_\theta \to
\infty$.
Для определенности будем считать, что $\theta\in
(0,1]$. При $\theta>1$ доказательство не меняется,
нужно только, согласно Теореме 1.3, число $\theta/2$
заменить на  $\theta -1$.
 Выделим из последовательности $\bold y^n$ слабо
сходящуюся  подпоследовательность в пространстве
$l^\theta_B$. Считаем что сама последовательность
слабо сходится  к некоторому элементу $\bold y\in
l^{\,\theta}_B$. Из слабой сходимости следует
покоординатная сходимость. Тогда из определения
множества $\Omega^{\,\theta}(h,r)$ и его замкнутости
следует, что $\bold y \in
\Omega^{\,\theta}(h,r)$. В силу Теоремы 2.1
найдется функция $\sigma\in
\Gamma^\theta$, такая, что $F\sigma =\bold y.$
Из слабой сходимости $\bold y^n\rightharpoonup\bold
y$ в $l^{\,\theta}_B$ следует сильная сходимость
$\bold y^n\to \bold y$  в норме $l^{\, \theta/2}_B$,
а из аналитичности (достаточно непрерывности)
отображения $F^{-1}:\,
\Omega^{\,\theta/2} \to \Gamma^{\,\theta/2} $
следует, что  $\|\sigma_n-\sigma\|_{\theta/2} \to 0$.
В силу Теоремы 1.3 имеем $\|\Phi\sigma_n\|_\theta
\leqslant\|\sigma_n\|_{\theta/2} \leqslant C$. Но тогда
(опять  используем Теорему 1.3 и свойство
ограниченности слабо сходящейся последовательности)
имеем
$$
\|T\sigma_n\|_\theta \leqslant \|\Phi
\sigma_n\|_\theta +\|\bold y^n\|_\theta \leqslant C+C
= 2C.
$$
Поскольку оператор $T: W^\theta_2 \to l^{\,\theta}_B$
есть изоморфизм, то $\|\sigma_n\|_\theta \leqslant
2C$. Это противоречие завершает доказательство
первого утверждения леммы.

Очевидно, что $\Psi =-T^{-1}\Phi F^{-1}$.
Следовательно, отображение $\Psi:
\Omega^{\,\theta} \to W^\tau_2$ является аналитическим как
композиция аналитических отображений. Из первого
утверждения леммы и равномерной ограниченности в
каждом шаре отображения $\Phi : \Omega^{\,\theta} \to
W^\tau_2$ получаем оценку $\|\Psi\bold
y\|_\tau\leqslant C$ для всех $\bold y\in
\Omega^{\,\theta}(r,h)$. Из равенства $\Psi(0) =0$  и аналитичности
отображения $\Psi$ получаем   оценку
\eqref{Psi}.
 Лемма доказана.\quad$\Box$

Следующее утверждение является совсем простым, но нам
удобно его отдельно сформулировать.

{\bf Лемма 2.9.} {\sl Пусть $X$, $X_1$ --- метрические
пространства, $X$ полно и функция $\Phi: X\to X_1$
 непрерывна на $X$. Если множество $U \subset X$
 предкомпактно в $X$, то $\Phi :U\to X_1$  равномерно
 непрерывна и равномерно ограничена.
}

{\bf Доказательство.} В условиях леммы замыкание $\overline {U}$ является компактом в $X$, а функция $\Phi: \overline
{U} \to X_1$ непрерывна. Поэтому утверждение следует из свойств непрерывных функций на компактах.\quad$\Box$

{\bf Лемма 2.10.} {\sl Пусть $\theta>0$. При любом $R>0$  справедлива
оценка
\begin{equation}\label{direct}
\|F'(\sigma)\|_\theta \leqslant C, \qquad \text{для всех}\ \ \sigma \in
\mathcal B_\Gamma^{\,\theta} (R),
\end{equation}
где постоянная $C$ зависит от $R$,  но не зависит от
$\sigma$. }

{\bf Доказательство.} Не ограничивая общности,
считаем, что $\theta\in (0,1]$. Если $\theta >1$, то
далее число $\theta/2$ нужно заменять на $\theta -1$.
Поскольку $F'=\Phi'+T$, достаточно доказать оценку
\eqref{direct},  в которой вместо $F$  участвует $\Phi$.
Согласно Теореме 2.3, отображение $\Phi :
W^{\theta/2}
\to l^{\,\theta}_B $  аналитично на замкнутом множестве
$\mathcal B_\Gamma^{\,\theta/2} (R_1)$ при любом
$R_1>0$, а потому числовая функция
$\|\Phi'(\sigma)\|_\theta$ непрерывна на этом
множестве. Из непрерывности вложения
$W^\theta_2\hookrightarrow W^{\theta/2}_2$ следует,
что найдется число $R_1= R_1(R, \theta)$,  такое, что
$\mathcal B_\Gamma^{\,\theta} (R)  \subset \mathcal
B_\Gamma^{\,\theta/2} (R_1)$. Здесь первое множество
компактно во втором, поэтому из Леммы 2.9 следует
оценка
\eqref{direct}, в которой $F$
надо заменить на $\Phi$.  Лемма доказана.\quad$\Box$

{\bf Лемма 2.11.} {\sl Пусть $\theta >0$. При любых
$r>0,\ h\in (0,\, 1/2)$  для обратного отображения
справедлива оценка

\begin{equation}\label{inverse}
\|(F^{-1})'(\bold y )\| \leqslant C, \qquad \text{для всех}
\ \, \bold y \in \Omega^\theta(r,h),
\end{equation}
где постоянная $C$ зависит от $r$ и $h$,  но не
зависит от $\bold y$. }

{\bf Доказательство.}
 Для определенности рассматриваем
случай $\theta\in (0,1]$. Рассуждаем аналогично.
Фиксируем числа $r>0,  \ \, h\in (0,\, 1/2)$.
Используя непрерывность вложения
$l^{\,\theta}_B\hookrightarrow l^{\, \theta/2}_B$
найдем число $r_1$,  такое,  что
$\Omega^{\,\theta}(r,h)
\subset
\Omega^{\,\theta/2}(r_1,h)$.
 Согласно
Лемме~2.8, отображение
$$
\Psi =-F^{-1} \Phi T^{-1}: \Omega^{\,\theta/2}(r_1,h)\to
W^\theta_2
$$
аналитично.  Поэтому числовая функция $\|\Psi'(\bold
y)\|_\theta $  непрерывна при $\bold y\in
\Omega^{\,\theta/2}(r_1,h)$. Воспользовавшись Леммой 2.9
и  компактностью вложения $\Omega^{\,\theta}(r,h)
\subset \Omega^{\,\theta/2}(r_1,h)$, получим оценку
\eqref{inverse}
в которой $F^{-1}$ заменено на $\Psi$. Поскольку $F^{-1} = T^{-1} + \Psi$, оценка сохраняется для $F^{-1}$. Лемма
доказана.\quad$\Box$

Теперь мы можем доказать основной результат этого
параграфа.

{\bf Теорема~2.12.} {\sl  Фиксируем $\theta >0$.
 Пусть последовательности
$\bold y, \bold y_1$  регуляризованных спектральных
данных лежат в  $\Omega_B^{\,\theta}(r,h)$. Тогда
прообразы $\sigma = F^{-1}_B\bold y,\ \,
\sigma_1 = F^{-1}_B\bold y_1$  лежат в
$\mathcal B_\Gamma^{\,\theta} (R)$ и справедливы
оценки
\begin{equation}\label{a}
C_1\|\bold y -\bold y_1\|_\theta \leqslant
\|\sigma -\sigma_1\|_\theta
\leqslant C_2\|\bold y -\bold y_1\|_\theta,
\end{equation}
где число $R$ и постоянные $C_1, C_2$ зависят только
$r$ и $h$. Число $R$   и постоянные $C_2, C^{-1}_1$
увеличиваются  при $r\to\infty$  или $h\to 0$.
Обратно, если $\sigma,
\sigma_1$ лежат в шаре $\mathcal B_\mathbb
R^{\,\theta}(R)$, то последовательности  $\bold y,
\bold y_1$
 регуляризованных спектральных данных этих функций
 лежат в  $\Omega^\theta(r,h)$
  и справедливы оценки

\begin{equation}\label{b}
C_1 \|\sigma -\sigma_1\|_\theta \leqslant\|\bold y
-\bold y_1\|_\theta \leqslant C_2 \|\sigma
-\sigma_1\|_\theta .
\end{equation}
Здесь числа $r>0, \ \, h\in (0, 1/2)$  и постоянные
$C_1$ и $C_2$ зависят только от $R$. Числа $ r,
\, h^{-1}, C_2 $ и $C^{-1}_1$  увеличиваются при
$R\to\infty$. }

{\bf Доказательство.}   Заметим, что множество
$\Omega^{\,\theta}(r,h) $ выпукло. Для
дифференцируемых функций на выпуклых множествах
справедлив аналог теоремы Лагранжа (см., например,
\cite[Следствие 12.2.8 гл. 12]{BS})
$$
\|\sigma -\sigma_1\|
\leqslant\sup\limits_{0<t<1}\|(F^{-1})'(t\bold y+(1-t)\bold
y_1)\|\ \cdot\,\|\bold y -\bold y_1\|.
$$
Тогда Лемма 2.11 влечет за собой оценку сверху в
неравенстве
\eqref{a}. Оценки сверху в  \eqref{b} получаются
аналогично из леммы 2.10. Оценки снизу в \eqref{a} и
\eqref{b} теперь следуют из оценок сверху и лемм 2.7 и 2.8.
Теорема доказана.\quad$\Box$

 Множества  $\mathcal B_\Gamma^\theta(R)$ в Теореме
2.12 можно заменить обычными шарами $\mathcal
B^\theta_{\mathbb R}(R)$, но тогда регуляризованные
спектральные данные нужно определить формулой
\eqref{shift},  где постоянная $c$  такова, что для
всех $\sigma \in \mathcal B^{\,\theta}_{\mathbb R}(R)
$ выполнена оценка $c\geqslant -
\mu_1(\sigma)-1/4$. В силу теоремы 3.1 такая
постоянная, зависящая только от $R$, существует. Это
замечание  вытекает из того, что при добавлении к
$\sigma$    функции $c(x-\pi)$ спектры операторов
$L_D$ и $L_{DN}$ сдвигаются на $c$, а разность
функций $\sigma ,\, \sigma_1 \in
\mathcal B^{\,\theta}_{\mathbb R}(R)$ совпадает с разностью
функций $\sigma +c(x-\pi), \ \sigma_1 +c(x-\pi)\in
\mathcal B_\Gamma^{\,\theta} (R)$.

\bigskip

\section{ Задача восстановления оператора $L_D$ по его
 спектральной функции.
Характеризация спектральных данных и равномерные
априорные оценки.}

\bigskip

 Общая схема доказательства аналогичных результатов
 для   задачи восстановления оператора $L_D$ по его
 спектральной функции
  остается прежней, хотя
 доказательства схожих по формулировке лемм
 проводятся по другому.
  В ходе изложения
мы сформулируем две леммы (Леммы 3.1 и 3.6),
доказательство которых носит технический характер. В
виду ограничения объема статьи мы укажем только путь,
на котором получаются доказательства, а детали и
подробные выкладки читатель может найти в нашей
электронной публикации \cite{SS6}.

Далее удобнее работать не с пространством $W^\theta_2\ominus\{1\}$, а с фактор-пространством $W^\theta_2 / \{1\}$,
считая, что все функции из $W^\theta_2$ определены с точностью до константы. Подразумеваем, что скалярное произведение
функций $f,g\in  W^\theta_2 / \{1\}$ определено равенством $(f,g)_\theta =(f_0, g_0)_\theta$,  где $f_0,g_0\in
W^\theta_2\ominus\{1\}$.
 Обозначим
через $\Gamma^{\,\theta}_D$ множество вещественных
функций
 $\sigma\in W^\theta_2 /\{1\}$, для
которых $\la_1(\sigma)\geqslant  1/2$, а через
$\mathcal B^{\,\theta}_\Gamma (R)$ --- пересечение
множества $\Gamma^{\,\theta}_D$ с замкнутым шаром
$\mathcal B^{\,\theta}_{\mathbb R}(R)$. Если
$\sigma\in\Gamma^{\,\theta}_D$,  то собственные
значения оператора $L_D$ подчинены условиям $1/2
\leqslant \lambda_1 <\la_2<\dots$.
Для регуляризованных спектральных данных эти
неравенства эквивалентны следующим:
\begin{equation}\label{s1}
s_2 \geqslant 0, \qquad s_{2k}-s_{2k+2}<1, \qquad
k=1,2,\dots .
\end{equation}
Условия неотрицательности всех нормировочных чисел
эквивалентны условиям
\begin{equation}\label{s2}
s_{2k-1}>-\pi/2,\qquad k=1,2,\dots.
\end{equation}
Последовательность $\{s_k\}_1^\infty$ принадлежит $l_2$, поэтому для любой вещественной функции $\sigma\in
\Gamma^{\,\theta}_D$ найдется число $h=h(\sigma)>0$,
такое, что
\begin{equation}\label{sh}
s_2 \geqslant 0, \qquad s_{2k}-s_{2k+2}\leqslant
1-h,\qquad s_{2k-1}\geqslant -\pi/2+h,
\qquad k=1,2,\dots .
\end{equation}

Фиксируем произвольные числа $r>0$ и $h \in (0,1)$.
Обозначим через $\Omega_D^\theta (r,h)$ совокупность
вещественных последовательностей $\{s_k\}_1^\infty$,
для которых выполнены неравенства \eqref{sh} и
которые лежат в замкнутом шаре радиуса $r$
пространства $l_D^{\,\theta}$,  т.е.
$\|\{s_k\}\|_\theta
\leqslant r$.
Через $\Omega_D^{\,\theta}$ обозначим множество  всех
вещественных последовательностей $\{s_k\}_1^\infty
\in l_D^{\,\theta}$, для которых справедливы неравенства
\eqref{s1} и \eqref{s2}.
Далее мы работаем только с отображением $F_D$ и, где
удобно, {\it будем опускать индекс $D$}. Вместо
$\Gamma_D^{\,\theta}$, $\Omega_D^{\,\theta}$ и
 $\Omega_D^\theta (r,h)$ всегда будем писать
$\Gamma^{\,\theta}$,\ \, $\Omega^{\,\theta}$ и
$\Omega^\theta (r,h)$  соответственно.

Для доказательства аналогов Теорем 2.1 и 2.2 нам
понадобится следующий важный результат, который дает
явное описание  прообраза отображения $F_D$ при
изменении только одной из координат в пространстве
$l^{\,\theta}_D$. Похожие формулы для задачи
восстановления по одному спектру имеются в книге
\cite{PoTr}. Но доказательство нашего результата
проводится на другом пути.

{\bf Лемма 3.1} {\sl Пусть $\{\la_k\}$ и $\{\al_k\}$ ---
собственные значения и нормировочные числа оператора
$L_D$ с вещественной функцией $\sigma\in
W_2^\theta\in \Gamma^{\,\theta}$, $\theta\ge0$. Тогда
для любого фиксированного $n\ge1$ и для любого
$t\in(\la_{n-1}-\la_n,\la_{n+1}-\la_n)$ существует
функция $\sigma(x,t)\in W_2^\theta$, такая, что
соответствующий оператор $L_D= L_D(\sigma)$ имеет
спектр $\{\la_k+t\delta_{kn}\}_1^\infty$ (здесь
$\delta_{kn}$ --- символ Кронекера) и нормировочные
числа $\{\al_k\}$. Далее, для любого фиксированного
$n\ge1$ и для любого $t\in(-\al_n,+\infty)$
существует функция $\sigma(x,t)\in W_2^\theta$,
такая, что оператор $L_D$, построенный по этой
функции, имеет спектр $\{\la_k\}_1^\infty$ и
нормировочные числа
$\{\al_k+t\delta_{kn}\}_1^\infty$.}

{\bf Доказательство.}
 Потенциал $\sigma(x,t)$ можно выписать в явном виде. В первом
случае, когда меняется собственное значение $\la_n$,
а  нормировочные числа и все другие собственные
значения остаются неизменными, положим

\begin{equation}\label{varsigma_la}
\sigma_n(x,t)=\sigma(x)- 2\tfrac{d}{dx}\ln G(x,t),
\end{equation}
где
\begin{multline}\label{G_n}
G(x,t) =\left(1+\al_n^{-1}\int_0^x y^2(\xi,\la_n+t)\,
d\xi\right)\left(1-\al_n^{-1}\int_0^x
y^2(\xi,\la_n)\,d\xi\right)\\
+\left(\al_n^{-1}\int_0^xy(\xi,\la_n+t)y(\xi,\la_n)\,
d\xi\right)^2.
\end{multline}
Здесь $y(x,\la)$ --- решение уравнения
$-y''+\sigma'y=\la y$ с начальными условиями
$y(0,\la)=0$, $y^{[1]}(0,\la)=\sqrt{\la}$. Во втором
случае, когда меняется только одно нормировочное
число $\alpha_n$,  положим
\begin{equation}\label{varsigma_al}
\sigma_n(x,t)=\sigma(x)- 2\tfrac{d}{dx}\ln G(x,t),\quad \
\text{где}\quad
G(x,t)=1+((\al_n+t)^{-1}-\al_n^{-1})\int_0^xy^2(\xi,\la_n)\,d\xi.
\end{equation}
Выписанные формулы получаются, если написать
уравнение Гельфанда--Левитана--Марченко в том виде, в
котором оно получено для потенциалов--распределений
Гринивым и Микитюком  \cite{HM2}. Если  искать
решения этого уравнения,
 удовлетворяющие условиям леммы,
 в виде линейной комбинации
двух функций (ср.
\cite[стр. 49-50]{Le2}),
  то получается система двух линейных уравнений,
которая решается явно. Подробности можно найти в
\cite{SS6}.\quad$\Box$

{\bf Лемма 3.2.} {\sl При любых $ \theta
\geqslant 0$\ отображение \ $F_D: \Gamma^{\,\theta}
\to \Omega^{\,\theta}_D$\ сюръективно.}

{\bf Доказательство.}
Сначала докажем лемму при $\theta<1/2$, когда
пространство  $l_D^{\,\theta}$ совпадает с
$l_2^{\,\theta}$. Воспользуемся приемом из
\cite{PoTr}. Согласно Теоремам 1.1 и 1.2 производная
по Фреше отображения $F_D$ в точке $\sigma =0$
совпадает с оператором $T_D$, который является
изоморфизмом. Поэтому для
 любого достаточно малого числа
$\varepsilon >0$ найдется такое $\delta>0$, что образ
шара $\|\sigma\|_\theta<\delta$ при отображении $F_D$
накрывает шар $\|s\|_\theta<\varepsilon$. При
$\theta<1/2$ пространство $l_D^{\,\theta}$ совпадает
с пространством $l_2^{\,\theta}$.  Для данного $\bold
s=\{s_k\}\in\Omega^\theta$ рассмотрим
последовательность
$$
\bold s^n =\{0,0,\dots,0,s_n,s_{n+1},\dots\},
$$
выбрав число $n$ так, чтобы $\|\bold s^n
\|_\theta<\varepsilon$. Тогда найдется единственная
функция $\sigma_n\in W_2^\theta$, образ $F(\sigma_n)$
 которой совпадает с $\bold s^n$. Применив лемму 3.1 ($n-1$)
раз, построим функцию
$\sigma\in\Gamma^{\,\theta}\subset W_{2,\mathbb
R}^\theta$, для которой $F \sigma =
\bold s$. Это и означает, что образ отображения $F$
содержит $\Omega^\theta$. Теперь при $\theta
\geqslant 1/2$ доказательство завершается с помощью
приема, примененного в доказательстве Теоремы 2.1
 Лемма доказана.\quad$\Box$

{\bf Лемма 3.3.} {\sl При любых
$ \theta \geqslant0$\ отображение \ $F:\,
\Gamma^{\,\theta} \to
\Omega^{\,\theta}$\ инъективно. }

{\bf Доказательство}. Инъективность этого отображения
при $\theta =0$ (а тогда при всех $\theta\geqslant
0$) доказана в работе Гринива и Микитюка \cite{HM2}.
Отметим также, что инъективность  следует из
формулируемой ниже Леммы 3.6  (для доказательства
нужно повторить рассуждения из Леммы 6 работы авторов
\cite{SS3}). Лемма доказана.\quad$\Box$

Обозначим через  $\widehat\Omega^{\,\theta}$
множество последовательностей $\{s_k\}_{k=1}^\infty
\in l_D^{\,\theta}$, для которых числа $\lambda_k =(s_k +k)^2$
вещественны. Повторив рассуждения, проведенные перед
доказательством теоремы 2.2, из Лемм 3.2 и 3.3
получаем аналог Теоремы 2.2.

{\bf Теорема 3.4.} {\sl При любых
$ \theta \geqslant 0$\ отображение \ $F_D:\,
W^\theta_{2,\mathbb R} /\{1\} \to
\widehat\Omega^{\,\theta}$\ есть биекция.
В частности, числа $\{\lambda_k\}_1^\infty$
 и $\{\alpha_k\}_1^\infty$ являются собственными
 значениями и нормировочными числами оператора $
L_D$, порождаемого функцией $\sigma\in W^\theta_{2,\mathbb R}$, если и только если последовательность $\{\lambda_k\}$
строго монотонна, числа $\{\alpha_k\}$ положительны и $\{s_k\}_1^\infty
\in l^{\,\theta}_D$.}

Из сформулированного утверждения следует также, что
отображение $F_D:\, \Gamma^{\,\theta}\to
\Omega^{\,\theta}$ есть биекция.  Отметим, что при
натуральных $\theta =1,2,
\dots$ аналог Теоремы 3.4, сформулированный на другом
языке, имеется в книге Фрайлинга и Юрко \cite
{FrYu}.

Аналитичность и явный вид производной по Фреше дает
следующая теорема.

{\bf Теорема 3.5.} {\sl Пусть $\theta\geqslant 0$ и
$\sigma\in\Gamma^{\,\theta}$. Тогда найдется
комплексная окрестность $U\in W_2^\theta$ точки
$\sigma$, такая, что отображение $F: U
\to l_D^\theta$  является вещественно аналитическим.
  В этой окрестности отображение $\Phi_D =F_D-T_D: U \to l^{\,\tau}_D
$, где $\tau$ определено в Теореме 1.3, также
является вещественно аналитическим. Производная в
точке $\sigma\in U$ определяется равенством
\begin{equation} \label{F'D}
F'_D(\sigma) f  = \left\{ (\varphi_k(x),
\overline{f(x)})
\right\}_{k=1}^\infty,
\end{equation}
где
\begin{equation}\label{phi_k}
\varphi_{2k-1}(x) =
2\al_k\la_k\frac{d}{d\la}(z(x,\la)z'(x,\la))\vert_{\la
=\la_k},
\quad\varphi_{2k}(x)=-\frac{y'_k(x)y_k(x)}{\al_k\sqrt{\la_k}},
\qquad k=1,2, \dots .
\end{equation}
Здесь $f\in W^\theta_2$ --- функция, на которую
действует оператор $ F'_D(\sigma): W_2^\theta\to
l^{\,\theta}_D$, $y_{n}=y(x,\la_n)$ --- собственные
функции оператора $L_D$, нормированные условиями
$y^{[1]}(0,\la_n)=\sqrt{\la_n}$, а $z(x,\la)$ ---
решение уравнения $-y''+\sigma'(x)y=\la y$ с
начальным условием $z(\pi,\la)=0$, нормированное
условием $\int_0^\pi z^2(x,\la)dx=\frac1\la$.
Утверждение об аналитичности (обычной) сохраняется,
если условие $\sigma\in\Gamma^{\,\theta}$ заменить
условием $\sigma\in W^{\,\theta}_{2,\mathbb R}$ и
потребовать, чтобы нуль не был собственным значением
оператора $L_D$.}

{\bf Доказательство.} Локальная дифференцируемость
отображения $F_D$ доказана в \S \ 6 нашей работы
\cite{SS5}. В этой же работе приведены
явные формулы для производной по Фреше, но они менее удобны, нежели \eqref{phi_k}. Переход от старых формул к новым
требует некоторой технической работы, см. \cite{SS6}.\quad$\Box$

{\bf Лемма 3.6}. {\sl Система функций $\{\varphi_k\}_1^\infty$, определенная равенствами \eqref{phi_k}, является
базисом Рисса в пространстве $L_2(0,
\pi) /\{1\}$. Биортогональная к ней система имеет вид
\begin{equation}\label{psiD}
\psi_{2k-1}(x)=\frac 2{\al_k^2}y_k^2(x),
\qquad\psi_{2k}(x)=-\frac{2\sqrt{\la_k}}{\al_k}\frac{d}
{d\la}\left(y^2(x,\la)\right)\vert_{\la =\la_k} \quad k=1,2,\dots ,
\end{equation}}
 а потому также является базисом Рисса.

{\bf Доказательство}
соотношений $(\varphi_k(x),\psi_n(x))=\delta_{kn}$
 при $k\ne n$ проводится так же, как в Лемме 6 работы
авторов \cite{SS5}. При $k=n$ проверка равенств усложняется, см.  \cite{SS6}.\quad$\Box$

Доказательства следующих двух теорем получаются дословным повторением доказательств Теорем 2.6 и 2.11 соответственно.

{\bf Теорема 3.7} {\sl  Пусть $\theta \geqslant 0.$
Для каждой точки $\bold y_0
\in \Omega^\theta =F_D(\Gamma^\theta)$
 существует ее комплексная окрестность $U(\bold y_0)$,
в которой определено обратное отображение
$F^{-1}_D(\bold y)$ и в которой это отображение имеет
комплексную производную по Фреше. Эта производная
имеет вид
\begin{equation*}
\left( F_D^{-1}\right)'(\bold y)= (F'_D)^{-1}(\bold y) =
\sum^\infty_{k=1} s_k\psi_k(x), \qquad \bold y= (s_1, s_2, \dots
).
\end{equation*}
Здесь $\{\psi_k(x)\}_1^\infty$ --- система,
биортогональная  к системе \eqref{phi_k}.}

{\bf Теорема 3.8.} {\sl
Утверждение Теоремы 2.12 сохраняет силу, если
отображение $F=F_B$  и множество $\Omega^{\,\theta}
(r,h)=\Omega_B^{\,\theta} (r,h)$  в ее формулировке
заменить на $F_D$  и $\Omega_D^{\,\theta }(r,h)$
соответственно. }

Авторы благодарят проф. Р.~О.~Гринива за прочтение
рукописи работы и полезные замечания.


\bigskip

 \address{A.М.Савчук,
 МГУ имени М.В.Ломоносова,
механико-математический ф-т, Ленинские Горы, Москва,
119992.
\email{artem\_savchuk@mail.ru}}

\address{A.А.Шкаликов, МГУ имени М.В.Ломоносова,
механико-математический ф-т, Ленинские Горы, Москва,
119992.
\email{ashkalikov@yahoo.com}
}
\end{document}